\documentclass[12pt,a4paper,leqno]{article}
\usepackage{amssymb, amscd, amsthm, amsmath, amstext,mathtools}
\usepackage{enumitem}
\usepackage[all]{xy}

\usepackage[colorlinks, linkcolor=blue, anchorcolor=green, citecolor=blue]{hyperref}

\newtheorem{thm}{Theorem}[section]
\newtheorem{coro}[thm]{Corollary}
\newtheorem{lemma}[thm]{Lemma}
\newtheorem{prop}[thm]{Proposition}

\theoremstyle{definition}

\newtheorem{example}[thm]{Example}
\newtheorem{ques}[thm]{Question}
\newtheorem{remark}[thm]{Remark}

\numberwithin{equation}{thm}
\begin{document}
\title{\textbf{Universality Criterion Sets for Definite Quadratic Lattices over Polynomial Rings}}
\author{Yong HU and Jing LIU}
\date{}
\maketitle 

\begin{abstract}
Let $k$ be a field of characteristic $\neq 2$. Under suitable assumptions on the dimensions of anisotropic quadratic forms over $k$ and $k(x)$, we show that a definite integral quadratic form over the polynomial ring $k[x]$ is universal if the form represents all scalars from $k$ and all scalar multiples of $x$. The coefficient field $k$ in our main result can be a finite field, a non-dyadic higher local field, a global function field, or a one-variable function field over $\mathbb{C}$ or $\mathbb{C}(t)$. We have a finite universality criterion set over $k[x]$ for some coefficient fields $k$, including some cases where $k$ has infinitely many square classes (e.g., $k=\mathbb{C}(t)$).
\end{abstract}



\section{Introduction}

For a quadratic form $f=f(x_1,\cdots,x_n)$ in $n$ variables over a ring $R$ and an element $\alpha\in R$, we say that $f$ represents $\alpha$ (over $R$) if there exist $a_1,\cdots, a_n\in R$ such that $f(a_1,\cdots, a_n)=\alpha$.

A positive definite quadratic form over the integers is called universal if it represents all positive integers. For example, the four-square theorem of Lagrange states that the form $x^2+y^2+z^2+w^2$ is universal. The famous 15-theorem of Conway and Schneeberger \cite{Con00,Bha00} identifies the set $\{1,2,3,5,6,7,10,14,15\}$ as a \textbf{universality criterion set} for positive definite quadratic forms with integer Gram matrix over $\mathbb{Z}$, in the sense that any such form which represents all elements in that set is universal. A theorem of Bhargava and Hanke \cite{BH11}, known as the 290-theorem, extends the 15-theorem to cover the case of integral quadratic forms that may have half-integers in their Gram matrices. These two theorems have found a large number of applications, demonstrating that it is of both theoretical and computational importance to have convenient universality criterion sets at our disposal.

Now let us consider the polynomial ring $k[x]$ with coefficients in a field $k$. We say that a quadratic form over $k[x]$ is universal if it represents all elements of the ring $k[x]$. Thus, a universality criterion set for integral quadratic forms over $k[x]$ can be defined as a subset of $k[x]\setminus \{0\}$ such that any integral quadratic form representing all elements of that set is universal.

Gerstein \cite{Ger04} pioneered the study of universal quadratic forms over polynomial rings. He defined a quadratic form over $k[x]$ to be definite if the form remains anisotropic over the field $k(\!(x^{-1})\!)$ (which is the completion of the rational function field $k(x)$ with respect to the negative-degree valuation). As a function field analogue of the 15-theorem and 290-theorem,  Gerstein's 4-conjecture speculated that $\{1,\delta,x,\delta x\}$ is a universality criterion set for definite quadratic forms over  $k[x]$, when $k$ is a finite field of odd characteristic and $\delta$ denotes a non-square in $k$. This conjecture has been confirmed with different methods in \cite{KWX08} and \cite{CD05}, and also in the thesis \cite{Bur06}.

Motivated by the above results, the current work aims to address the problem of finding universality criterion sets for definite quadratic forms over a polynomial ring $k[x]$ for more general coefficient fields $k$. Under suitable assumptions on $k$, we show that it will still be sufficient to represent all scalars from $k$ and all scalar multiples of $x$ for a definite integral quadratic form over $k[x]$ to be universal, just as in the case with $k$ finite.

Similar to the cases already discussed above, a quadratic form over any field $k$ is called universal if it represents all elements of $k^*=k\setminus\{0\}$, and a universality criterion set for quadratic forms over $k$ is a subset of $k^*$ that is sufficient to test the universality of quadratic forms over $k$.

The $u$-invariant $u(k)$ of a field $k$ is defined as the supremum of the dimensions of anisotropic quadratic forms over $k$ (cf. \cite[Definition\;XI.6.1]{Lam05}). We say that  $k$ has UUD (short for unique universality dimension) if its $u$-invariant $u(k)$ is finite and every anisotropic universal quadratic form over $k$ is of dimension $u(k)$.

\begin{thm}\label{thm-main}
    Let $k$ be a field of characteristic $\neq 2$. Assume that $k$ has \textup{UUD} and that $u(k(x))=2u(k)$. Then for any universality criterion set $S$ for quadratic forms over $k$, the set
\[
S\cup Sx=S\cup\{ax\,|\,a\in S\}
\] is a universality criterion set for definite integral quadratic forms over the polynomial ring $k[x]$.
\end{thm}

Clearly, we can take $S$ to be any complete system of representatives of the group $k^*/k^{*2}$ of square classes. Therefore, Theorem~\ref{thm-main} gives a finite universality criterion set over $k[x]$ if  $k^*/k^{*2}$ is finite. In some examples with $k^*/k^{*2}$ infinite, we can still find a finite set $S$ as universality criterion set over $k$. See Remark\;\ref{4.5}.

As we will explain in Example\;\ref{4.3}, typical examples of fields $k$ satisfying the hypotheses in Theorem~\ref{thm-main} include finite fields, non-archimedean local fields, global function fields, non-dyadic higher local fields, and one-variable function fields over $\mathbb{C}$, $\mathbb{C}(t)$ or  $\mathbb{C}(\!(t)\!)$. In particular, the present work gives yet another proof of the aforementioned 4-conjecture of Gerstein, previously proved in \cite{CD05,Bur06,KWX08}. Also, for any odd prime $p$, with $\delta$ denoting a non-square unit in the field $\mathbb{Q}_p$ of $p$-adic numbers, one obtains $\{1,\delta,p,\delta p,x,\delta x,px,\delta px\}$ as a universality criterion set for $\mathbb{Q}_p[x]$. This might analogously be referred to as the 8-theorem.

\medskip

Let us sketch the strategy of proof of our main theorem and outline the organization of the paper as follows: In \S~\ref{sec-suc. min.} we extend the theory of successive minima for definite quadratic forms over $k[x]$, first developed in \cite{Ger03} for $k$ finite, and now in full generality. The theory is then combined with the UUD property in \S~\ref{sec-UUD}  to deduce a simple reduced form for any integral quadratic form representing monomials of degree $\leqslant 1$. Finally in \S~\ref{sec-main}, we conclude the proof by using the generalized Cassels--Pfister theorem established in \cite{Cla12}, a key notion involved here being the so-called ADC property of integral quadratic forms.

\medskip

\noindent {\bf Notation and Conventions}. Throughout what follows, let $k$ be a field of characteristic different from $2$.

We will use the geometric language of quadratic spaces and lattices, standard notation and terminology about which will be adopted from \cite{Lam05} and \cite{OMeara00}.
In particular, a quadratic space over  $k$ is a pair $(V,\varphi)$ consisting of a finite dimensional $k$-vector space $V$ and a quadratic form $\varphi\colon V\to k$. All the quadratic spaces considered in this paper are assumed to be non-degenerate unless otherwise stated.

A scalar $a\in k$ is said to be represented by $(V,\varphi)$ if there exists a nonzero vector $v\in V$ with $\varphi(v)=a$. A quadratic space over $k$ is called universal if it represents every element in $k^*$. On the other hand, a quadratic space is called isotropic over $k$ if it represents zero, and anisotropic otherwise.

The natural base change of a quadratic space $(V,\varphi)$ over $k$ along a field extension $k'/k$ is given by $(V_{k'}=V\otimes_{k}k',\varphi')$ with $\varphi'(v\otimes a)\coloneqq a^2 \varphi(v)$ for all $v\in V$ and $a\in k'$. This amounts to viewing the same quadratic form $\varphi$ as a quadratic form $\varphi_{k'}$ with coefficients in $k'$.

When there is no risk of confusion, we often use the terms \emph{quadratic form} and \emph{homogeneous polynomial of degree $2$} synonymously.
We write $\langle a_1,\ldots,a_n\rangle$ for the diagonal quadratic form $\sum_{i=1}^{n}a_i x_i^2$.

\section{Successive Minima for Definite Quadratic Lattices}\label{sec-suc. min.}

In \cite[\S\,3, Thm.\;1]{Ger03}, a theory of successive minima for definite quadratic lattices over the polynomial ring $k[x]$ is established when $k$ is a finite field. Our goal here is to show that this theory of successive minima is valid for arbitrary $k$.

\begin{lemma}\label{2.1}
Let $(K,v)$ be a henselian valuation field with valuation ring $\mathcal{O}_K$ and residue field $k$ (of characteristic $\neq 2$).
For an element $u\in\mathcal{O}_K$, let $\bar u$ denote its canonical image in the residue field $k$.

  \begin{enumerate}[label={{\upshape(\arabic*)}}]
    \item A unit $u$ in $\mathcal{O}_K$ is a square in $K$ if and only if $\bar u$ is a square in $k$.
    \item Let $u_1,\cdots, u_n$ be units in the valuation ring $\mathcal{O}_K$.

   Then the diagonal quadratic form $\langle u_1,\cdots, u_n\rangle$ is isotropic over $K$ if and only if the form $\langle \bar u_1,\cdots,\bar u_n\rangle$ is isotropic over $k$.
  \end{enumerate}
\end{lemma}
\begin{proof}
   This is a well known consequence of Hensel's lemma. In this paper we only use the case where $(K,\,v)$ is a complete discrete valuation field. A proof in that special case  can be found in \cite[p.146, Lemma\;VI.1.1 and Prop.\;VI.1.6]{Lam05}.
\end{proof}
%

\begin{lemma}\label{2.2}
  Let $(K,v)$ be a henselian valuation field with residue field $k$. Let $\varphi$ be a quadratic form of dimension $n\geqslant 1$ over $K$. Suppose that the Gram matrix $A=(a_{ij})\in \mathrm{M}_n(K)$ of $\varphi$ in some basis satisfies $v(a_{ii})<v(a_{ij})$ for all $1\leqslant i\neq j\leqslant n$.

  Then $\det(\varphi)=\prod^n_{i=1}a_{ii} \in K^*/K^{*2}$ and $\varphi\cong \langle a_{11},\cdots, a_{nn}\rangle$.
\end{lemma}
\begin{proof}
  This is a restatement of \cite[\S\,2, Lemma\;2]{Ger03}, whose proof is still valid without assuming $k$ finite. Indeed, using Lemma\;\ref{2.1} (1), we can prove $\det(\varphi)=\prod^n_{i=1}a_{ii}$ in the same way as in Gerstein's paper. Let $\psi$ be the quadratic form defined by the $(n-1)\times (n-1)$ submatrix in the upper left corner of $A$. Then $\det(\psi)=\prod_{i=1}^{n-1} a_{ii}$ and in particular $\psi$ is non-degenerate, whence $\varphi\cong \psi\perp\langle d\rangle$ for some $d\in K^*$. By induction on $n$, we have $\psi\cong \langle a_{11},\cdots, a_{n-1,\,n-1}\rangle$. Comparing the determinants, we get $d=a_{nn}\in K^*/K^{*2}$; hence $\varphi\cong \psi\perp\langle d\rangle\cong \langle a_{11},\cdots, a_{nn}\rangle$, as desired.
\end{proof}

 For all $f,\,g\in k[x]$ with $g\neq 0$, we define $\deg(f/g):=\deg(f)-\deg(g)$. Then the map $\alpha\mapsto v_{\infty}(\alpha):=-\deg(\alpha)$ is a discrete valuation on the field $F:=k(x)$. We denote by $F_{\infty}$ the completion of $F$ at $v_{\infty}$.

\begin{lemma}\label{2.3}
  Let $\varphi$ be a quadratic form over the field $F=k(x)$. Suppose that the Gram matrix $A=(a_{ij})\in \mathrm{M}_n(F)$ of $\varphi$ in some basis has all its entries in $k[x]$ and satisfies $\deg(a_{ii})>\deg(a_{ij})$ for all $1\leqslant i\neq j\leqslant n$. Let $\lambda_ix^{\nu_i}$ be the leading term of $a_{ii}$, with $\lambda_i\in k$ and $\nu_i\in\mathbb{N}$.

  Then $\varphi_{F_\infty}\cong \langle a_{11},\cdots, a_{nn}\rangle\cong \langle \lambda_1x^{\nu_1},\cdots,\lambda_nx^{\nu_n}\rangle$ over the completion $F_{\infty}$.
\end{lemma}
\begin{proof}
  As in the proof of \cite[\S\,2, Lemma\;3]{Ger03}, this follows from Lemmas\;\ref{2.1} and \ref{2.2}.
\end{proof}

Let us recall some definitions for quadratic lattices over the polynomial ring $R:=k[x]$. By a (quadratic) lattice over $R$ we mean a finitely generated free $R$-module $L$ together with a non-degenerate quadratic form $\varphi$ defined on the vector space $FL=F\otimes_RL$. Put $V=FL$. We will also say that $L$ is an $R$-lattice on the quadratic space $(V,\,\varphi)$. The dimension $\dim (L)$ of $L$ is defined as the dimension of the $F$-vector space $V=FL$. If $\varphi(L)\subseteq R$, we say that the $R$-lattice $L$ is \emph{\textbf{integral}}.

For any $R$-lattice $L$ with associated quadratic form $\varphi$, the \emph{\textbf{minimum}} $\min(L)$ of $L$ is defined by
\[
\min(L):=\min\{\deg(\varphi(v))\,|\,v\in L\,,\; \varphi(v)\neq 0\}\,.
\]A \emph{\textbf{minimum vector}} of $L$ is a vector $v\in L$ such that $\deg(\varphi(v))=\min(L)$.

By a theorem of Djokovi\'c (cf. \cite[\S\,4, Thm.\;1]{Djo76} or \cite[\S\,1]{Ger03}), every $R$-lattice on an anisotropic quadratic $F$-space has a \emph{\textbf{reduced basis}}, i.e., an $R$-basis $v_1,\cdots, v_n$ in which the Gram matrix $(a_{ij})\in \mathrm{M}_n(F)$ satisfies the following two conditions:

\begin{enumerate}[label={{\upshape(\arabic*)}}]
\item
Diagonal non-decreasing: $\deg(a_{11})\leqslant \cdots\leqslant \deg(a_{nn})$.

\item
Diagonal domination: $\deg(a_{ii})>\deg(a_{ij})$ whenever $i\neq j$.
\end{enumerate}

A quadratic space $V$ over the field $F=k(x)$ is called \emph{\textbf{definite}} if $V_{\infty}:=V\otimes_FF_{\infty}$  is anisotropic. A quadratic $R$-lattice is called \emph{\textbf{definite}} if the quadratic $F$-space it spans is definite.

\medskip

When $k$ is finite, the following theorem is proved  in \cite[\S\,3, Thm.\;1]{Ger03}.

\begin{thm}\label{2.4}
Let $L$ be an integral lattice over $R=k[x]$ on a definite quadratic space with associated quadratic form $\varphi$, and let $A=(a_{ij})\in\mathrm{M}_n(R)$ be the Gram matrix of $L$ in some reduced basis $v_1,\cdots, v_n$. Let $v\in L$ be a nonzero vector and write $v=\sum^n_{i=1}\alpha_iv_i$, with $\alpha_i\in R$.

\begin{enumerate}[label={{\upshape(\arabic*)}}]
  \item  The leading term of the polynomial $\varphi(v)\in R$ is equal to the leading term of  $\sum^n_{i=1}\alpha_i^2 a_{ii}$, and $\deg(\varphi(v))=\max\limits_{1\leqslant i\leqslant n}\deg(\alpha_i^2a_{ii})$.
  \item $\deg(a_{11})=\min(L)$, and for $i>1$,
  \[
  \deg(a_{ii})=\min\{\deg(\varphi(u))\,|\, u\in L\,,\,u\notin Rv_1+\cdots+Rv_{i-1}\}\,.
  \]
  \item Let $1\leqslant r\leqslant n$ be the index such that
  \[
  \min(L)=\deg(a_{11})=\cdots=\deg(a_{rr})<\deg(a_{r+1,\,r+1})\,.
  \] Then a nonzero vector $v\in L$ is a minimum vector if and only if   $v\in kv_1+\cdots+kv_r$.
  \item Let $\{\deg(a_{11}),\cdots, \deg(a_{nn})\}=\{d_1,\cdots, d_t\}$ with $d_1<\cdots<d_t$, and put $d_0=-1$ and $d_{t+1}=+\infty$. For each $0\leqslant j\leqslant t$, let $L_j$ be the sublattice of $L$ generated by the set   $\{u\in L\,|\,\deg(\varphi(u))\leqslant d_j\}$.

      Then $L_0=0\subseteq L_1\subseteq \cdots\subseteq L_t=L$, and for $0\leqslant j\leqslant t$, the sublattice $L_j$ is generated by the set      $\{v_i\,|\, \deg(a_{ii})<d_{j+1}\}$.

      Consequently,
      \[
\#\{1\leqslant i\leqslant n\,|\, \deg(a_{ii})=d_j\}=\dim L_{j}-\dim L_{j-1}\\      \text{ for all } 1\leqslant j\leqslant t\,.
      \]
  \item   The sequence $(\deg(a_{11}),\cdots, \deg(a_{nn}))$ is independent of the choice of reduced basis for $L$. It is called the sequence of \emph{\textbf{successive minima}} of $L$.

\end{enumerate}
\end{thm}
\begin{proof}
As in the case treated in \cite[\S\,3, Thm.\;1]{Ger03},  (2)--(5) can be deduced from (1), and putting $M=\max\limits_{1\leqslant i\leqslant n}\deg(\alpha_i^2a_{ii})$, we have
$\deg(\alpha_i\alpha_ja_{ij})<M$ for all $1\leqslant i\neq j\leqslant n$, by the diagonal domination property of the matrix $A$.

For each $i$, let $\varepsilon_i\in k^*$ (resp. $\lambda_i\in k^*$) be the leading coefficient of $\alpha_i$ (resp. $a_{ii}$).  Let $J\subseteq\{1,\cdots, n\}$ be the set of indices $i$ such that $\deg(\alpha_i^2a_{ii})=M$. Then
\[
\varphi(v)=\sum^n_{i=1}\alpha_i^2a_{ii}+2\sum_{1\leqslant i<j\leqslant n}\alpha_i\alpha_ja_{ij}=\bigg(\sum_{i\in J}\varepsilon_i^2\lambda_i\bigg)x^M+ r_1
\]and
\[
\sum^n_{i=1}\alpha_i^2a_{ii}=\bigg(\sum_{i\in J}\varepsilon_i^2\lambda_i\bigg)x^M+ r_2\,,
\]where $r_1$ and $r_2$ are polynomials  of degree $<M$. Therefore, to finish the proof, it remains to show that $\sum_{i\in J}\varepsilon_i^2\lambda_i\neq 0$.

Indeed, if $\sum_{i\in J}\varepsilon_i^2\lambda_i=0$, then the $|J|$-dimensional diagonal quadratic form $\psi:=\langle \lambda_i\rangle_{i\in J}$ with $\lambda_i,\,i\in J$ as its diagonal entries is isotropic over $k$ and hence also over $F=k(x)$. Since $\deg(a_{ii})\equiv M\pmod{2}$ for all $i\in J$, the diagonal form $\psi':=\langle \lambda_ix^{\deg(a_{ii})}\rangle_{i\in J}$ is also isotropic over $F$.  But Lemma\;\ref{2.3} shows that $\psi'_{F_{\infty}}$ is a subform of $\varphi_{F_{\infty}}$. So $\varphi_{F_\infty}$ is isotropic. This contradicts the definiteness assumption on $\varphi$.
\end{proof}

\section{Consequences of the UUD Property}\label{sec-UUD}

Recall that the $u$-\emph{\textbf{invariant}} $u(k)$ of  $k$ is defined as the supremum of the dimensions of anisotropic quadratic spaces over $k$ (cf. \cite[Definition\;XI.6.1]{Lam05}). It is also the minimum of all integers $n\geqslant 1$ such that all $n$-dimensional quadratic spaces over $k$ are universal.

We say that the field $k$ has the \emph{\textbf{Unique Universality Dimension}} (UUD) property, or simply that $k$ has UUD, if $u(k)$ is finite and if the dimension of every anisotropic universal quadratic space over $k$  equals $u(k)$.

\begin{example}\label{3.1}
 If $u(k)\leqslant 4$, then $k$ has UUD.

 Indeed, we have $u(k)\in \{1, 2,\,4\}$ (\cite[Prop.\;XI.6.8]{Lam05}).  The result is clear if $u(k)\leqslant 2$. If $u(k)=4$, this follows from the well known fact that a ternary form of the type $\langle 1,\,a,\,b\rangle$ is anisotropic if and only if the quaternary form $\langle 1,\,a,\,b,\,ab\rangle$ is anisotropic.

 We see therefore that finite fields, non-archimedean local fields, global fields without real places, the rational function fields $\mathbb{C}(t), \mathbb{C}(t_1,t_2),\,\mathbb{C}(\!(t_1)\!)(t_2)$, and the Laurent series fields $\mathbb{C}(\!(t)\!),\,\mathbb{C}(\!(t_1)\!)(\!(t_2)\!),\,\mathbb{C}(t_1)(\!(t_2)\!)$ and $\mathbb{C}(\!(t_1,t_2)\!)$ all have the UUD property.
\end{example}

\begin{prop}\label{3.2}
Let $K$ be a complete discrete valuation field with residue field $k$ (of characteristic $\neq 2$).
Then $K$ has \textup{UUD} if and only if $k$ has \textup{UUD}.
\end{prop}
\begin{proof}
 We have $u(K)=2u(k)$ (cf. \cite[Example\;XI.6.2\,(7)]{Lam05}), so we may assume that $u(k)$ is finite.

 Let $\bar\varphi$ be an anisotropic universal quadratic form of dimension $d$ over $k$. Choose a diagonalization $\bar\varphi=\langle \bar a_1,\cdots \bar a_d\rangle$ and lift it to a form $\varphi=\langle a_1,\cdots,a_d\rangle$, where each $a_i$ is a unit in the valuation ring $\mathcal{O}_K$ that lifts $\bar a_i$. Then,  by Springer's theorem (cf. \cite[Prop.\;VI.1.9]{Lam05}), the form $f:=\varphi\perp \pi\varphi$, where $\pi$ is a uniformizer of $K$,  is anisotropic over $K$. On the other hand, for all units $u$ in $\mathcal{O}_K$, $\varphi$ represents $u$ and $\pi\varphi$ represents $u\pi$ by Lemma\;\ref{2.1}. It follows that $f$ is universal over $K$. Thus, if $K$ has UUD, then we must have $2d=u(K)=2u(k)$, i.e., $d=u(k)$. This shows that $k$ has UUD if $K$ does.

 Conversely, let $q$ be an anisotropic universal form over $K$. We can assume that $q=\psi\perp \pi\phi$, where $\psi$ and $\phi$ are diagonal forms whose diagonal entries are all units in $\mathcal{O}_K$. By Springer's theorem, the forms $\psi$ and $\phi$ are both anisotropic over $K$, and hence the residual forms $\bar{\psi}$ and $\bar{\phi}$ are anisotropic over $k$ by Lemma~\ref{2.1}. Moreover, for any $\bar{\delta}\in k^*$, letting $\delta\in \mathcal{O}_K$ be a unit lifting $\bar{\delta}$, another application of Springer's theorem together with Lemma~\ref{2.1} shows that $q$ represents $\delta$ over $K$ if and only if $\bar{\psi}$ represents $\bar{\delta}$ over $k$; and $q$ represents $ \pi\delta$ if and only if $\bar{\phi}$ represents $\bar{\delta}$ over $k$.  Since $q$ is universal over $K$, it follows that both $\bar{\psi}$ and $\bar{\phi}$ are universal over $k$. So, if $k$ has UUD, then $\dim\psi=\dim\phi=\dim\bar{\psi}=\dim\bar{\phi}=u(k)$, whence $\dim q=2u(k)=u(K)$. This proves the sufficiency part of the proposition.
\end{proof}

\begin{remark}\label{3.3}
  If, as in \cite{GesquiereVanGeel92} and \cite{Cassady26JPAA}, we define the $m$-\emph{\textbf{invariant}} $m(k)$ by
  \[
  m(k)=\inf\{\dim(\varphi)\,|\,\varphi \text{ an anisotropic universal quadratic form over } k\}\in \mathbb{N}\cup\{+\infty\},
  \]then $k$ has UUD if and only if $m(k)=u(k)<+\infty$. So Proposition\;\ref{3.2} can also be deduced directly from \cite[Cor.\;2.8]{Cassady26JPAA}.

  According to \cite[Cor.\;2.4]{Cassady26JPAA}, if $k$ has UUD, then $u(k)$ is a power of $2$.
\end{remark}

\begin{lemma}\label{3.4}
 Suppose that $k$ is a field such that $u(k)=m<+\infty$. 

  \begin{enumerate}[label={{\upshape(\arabic*)}}]
    \item Every definite quadratic space over $k(x)$ has dimension $\leqslant 2m$.
    \item Let $S\subseteq k^*$ be a universality criterion set for quadratic forms over $k$. Suppose that $k$ has \textup{UUD}.

     Then a definite  quadratic space over $k(x)$ which represents $\delta$ and $\delta x$ for all $\delta\in S$ must be of dimension $2m$.
  \end{enumerate}
\end{lemma}
\begin{proof}
  The completion $F_{\infty}$ of $F=k(x)$ at the discrete valuation $v_{\infty}$ is a complete discrete valuation field with residue field $k$. We have $u(F_{\infty})=2u(k)=2m$. This implies (1).

  To prove (2), first note that any square class in $F_{\infty}=k(\!(x^{-1})\!)$ is represented by an element of the form $\delta$ or $\delta x$ for some $\delta\in k$. Thus, using Springer's theorem (\cite[Prop.\;VI.1,9]{Lam05}) one can easily show that the set
   \[
   S\cup Sx=S\cup\{\delta x\,|\,\delta\in S\}
   \]is a universality criterion set for quadratic forms over $F_{\infty}$. So, if $V$ is a definite  quadratic space over $k(x)$ that represents all elements in $S\cup Sx$, then $V_{\infty}$ is an anisotropic universal space over $F_{\infty}$. By Proposition\;\ref{3.2}, $F_{\infty}$ has UUD, so we must have $\dim V=u(F_\infty)=2m$.
\end{proof}

\begin{coro}\label{3.5}
Let $S\subseteq k^*$ be a universality criterion set for quadratic forms over $k$. Let $L$ be a definite integral quadratic lattice of dimension $n$ over $k[x]$ and let $(d_1,\cdots, d_n)$ be its sequence of successive minima. Set
\[
r:=\#\{1\leqslant i\leqslant n\,|\,d_i=0\} \quad\text{and}\quad s:=\#\{1\leqslant i\leqslant n\,|\, d_i=1\}\,.
\]

\begin{enumerate}[label={{\upshape(\arabic*)}}]
  \item If $L$ represents an element in $k^*$, then $1\leqslant r\leqslant u(k)$.
  \item If $k$ has \textup{UUD} and $L$ represents all the elements in $S$, then $r=u(k)$.
  \item If $L$ represents an element in $k^*x=\{\delta x\,|\,\delta\in k^*\}$, then $1\leqslant s\leqslant u(k)$.
  \item If $k$ has \textup{UUD} and $L$ represents all the elements in $Sx$, then $s=u(k)$.
  \item If $k$ has \textup{UUD} and $L$ represents all the elements in $S\cup Sx$, then $r=s=u(k)$ and $n=2u(k)$, in other words, the sequence of successive minima of $L$ is
  \[
  (d_1,\cdots, d_n)=(0,\cdots, 0,\,1,\cdots, 1)\,,
  \]where both $0$ and $1$ occur $u(k)$ times.
\end{enumerate}
\end{coro}
\begin{proof}
  Let $\varphi$ be the quadratic form associated to $L$. Let $v_1,\cdots, v_n$ be a reduced basis of $L$ and let $A=(a_{ij})\in \mathrm{M}_n(k[x])$ be the Gram matrix of $\varphi$ in this reduced basis.
\begin{enumerate}[label={{\upshape(\arabic*)}}]
\item
Suppose that $L$ represents some element in $k^*$. Then $\min(L)=0$, so $r\geqslant 1$ by Theorem\;\ref{2.4} (2).

Let $\varphi_0= \langle a_{11},\cdots,a_{rr} \rangle$ be the restriction of $\varphi$ to the $k$-vector space spanned by $v_1,\cdots, v_r$. This $r$-dimensional form is anisotropic over $k$ since $L$ is anisotropic.
So we get $r\leqslant u(k)$.

\item
By Theorem\;\ref{2.4} (3), an element in $k^*$  is represented by $L$ over $k[x]$ if and only if it is represented by $\varphi_0$ over $k$.
So the assumption implies that $\varphi_0$ is universal over $k$.
Hence, the UUD property of $k$ yields $r=u(k)$, as desired.

\item
Let $\psi$ be the quadratic form over $F=k(x)$ obtained by restricting $\varphi$ to the $F$-vector space spanned by $v_{r+1},\cdots,v_{r+s}$. For $r+1\leqslant i\leqslant r+s$, let $\lambda_{i}$ denote the leading coefficient of the degree $1$ polynomial $a_{ii}$,
and define the diagonal form $\phi:=\langle \lambda_{r+1},\cdots, \lambda_{r+s}\rangle$ over $k$.
Then $\psi_{F_\infty}\cong \langle \lambda_{r+1}x,\cdots, \lambda_{r+s}x\rangle_{F_\infty}=x \phi_{F_\infty}$ by Lemma\;\ref{2.3}.
Since $\varphi_{F_\infty}$ is anisotropic ($L$ being definite), $\phi$ must be anisotropic over $k$.
Hence $s\leqslant u(k)$.

If $L$ represents $\delta x$ for some $\delta\in k^*$, then Theorem\;\ref{2.4} (1) shows that $\delta$ is represented by $\phi$ over $k$, hence $s\geqslant1$.

\item
If $L$ represents $\delta x$ for all $\delta\in S$, Theorem\;\ref{2.4} (1) shows that the form $\phi$ represents $S$, and hence $\phi$ is universal over $k$. Since $\phi$ is also anisotropic as we have noted, we conclude $s=u(k)$ if $k$ has UUD.

\item
This follows by combining (2), (4) and Lemma\;\ref{3.4} (2).
\qedhere
\end{enumerate}
\end{proof}


\section{Proof of Main Theorem, Examples and Remarks}\label{sec-main}

We can now give a proof of our main result, Theorem\;\ref{thm-main}.

Let us assume that the field $k$ (of characteristic $\neq 2$) has \textup{UUD} and that $u(k(x))=2u(k)$. Let $S\subseteq k^*$ be a universality criterion set for quadratic forms over $k$. Then Theorem\;\ref{thm-main} asserts that $S\cup Sx$ is a universality criterion set for definite integral quadratic lattices over the polynomial ring $k[x]$; that is, if a definite integral quadratic $k[x]$-lattice represents all elements in $S\cup Sx$, then the lattice represents all nonzero elements of $k[x]$.

\begin{proof}[Proof of Theorem\;$\ref{thm-main}$]
    Let $L$ be a definite integral quadratic $k[x]$-lattice that represents $S\cup Sx$. Then $\dim(L)=2u(k)$, and all the entries of the Gram matrix of $L$ with respect to a reduced basis have degree  $0$ or $1$, by Corollary~\ref{3.5} (5). The generalized Cassels--Pfister theorem \cite[Thm.~10]{Cla12} then tells us that (the integral quadratic form over $k[x]$ corresponding to) $L$ is an ADC form (cf. \cite[p.~144]{Cla12}), in the sense that it represents over the ring $k[x]$ any element of $k[x]$  that it represents over the field $k(x)$. But since $\dim(L)=u(k(x))$, the underlying quadratic space of $L$ is universal over $k(x)$; therefore $L$ is universal as a quadratic $k[x]$-lattice.
\end{proof}

\begin{prop}\label{4.1}
Let $k$ be a field (of characteristic $\neq 2$) such that for every finitely generated field extension $\ell/k$ of transcendence degree $d\leqslant 1$, one has $u(\ell)\leqslant 2^du(k)$. Let $K$ be a complete discrete valuation field with residue field $k$.

Then for every finitely generated field extension $L/K$ of transcendence degree $d\leqslant 1$, one has $u(L)\leqslant 2^du(K)$. In particular, $u(K(x))=2u(K)$.
\end{prop}
\begin{proof}
This is an easy consequence of \cite[Thm.\;4.10]{HHK09}.
\end{proof}

\begin{remark}\label{4.2}
	Let $k_0$ be a field of characteristic $\neq 2$ and let $k/k_0$ be a finitely generated field extension of transcendence degree $d$. Then $u(k)=2^d u(k_0)$ if $k_0$ is an
	algebraically closed field, a finite field, or a non-archimedean local field. In particular, this implies that $k$ satisfies the hypothesis of Proposition\;\ref{4.1}.
	
	Indeed, if $k_0$ is a non-archimedean local field of characteristic 0 (including the dyadic case), this is proved in \cite[Thm.\;3.4]{Leep13}.
	
	In the other cases, one can use the Tsen--Lang--Nagata theory of $C_i$-fields (see e.g., \cite[\S\,2.15]{Schar}) to show $u(k) \leqslant 2^d u(k_0)$. The reverse inequality is also well-known to the experts; it follows for example from \cite[Prop.\;5.2]{BGVG14}.
\end{remark}

\begin{example}\label{4.3}
By Example\;\ref{3.1}, Props.\;\ref{3.2} and \ref{4.1} and Remark\;\ref{4.2}, the field $k$ in Theorem\;\ref{thm-main} can be any of the following fields:
  \begin{enumerate}[label={{\upshape(\arabic*)}}]
    \item A finitely generated extension of transcendence degree $\leqslant 2$ over an algebraically closed field of characteristic $\neq 2$.
    \item A finitely generated extension of transcendence degree $\leqslant 1$ over $C(\!(t)\!)$, where $C$ is an algebraically closed field of characteristic $\neq 2$.
    \item A finite field or a global function field of characteristic $\neq 2$.
    \item A non-archimedean local field of characteristic $\neq 2$.
\item A complete discrete valuation field whose residue field is of the type as in the previous cases.
\end{enumerate}
\end{example}

Let $m\geqslant 1$ be an integer. As is customary, by an $m$-\emph{\textbf{local field}} we mean a field $\ell$ such that there is a sequence of fields $\ell_1,\ldots,\ell_m$ with $\ell_1$ a non-archimedean local field in the usual sense, $\ell_m=\ell$, and $\ell_i$ a complete discrete valuation field with residue field $\ell_{i-1}$ for $2\leqslant i\leqslant m$. By Example\;\ref{4.3} (5), in Theorem\;\ref{thm-main} we can take $k$ to be an $m$-local field $\ell$ with $\mathrm{char}(\ell_1)\neq 2$. In fact, by \cite[Cor.\;4.8 and Thm.\;4.10]{Cassady26JPAA}, we can also take $k=\ell(t)$ for an $m$-local field $\ell$ with $\mathrm{char}(\ell_1)\neq 2$.

 It would be interesting to know whether the following question has a positive answer:

\begin{ques}
 For $n\geqslant 3$, do the fields
 \[
 \mathbb{C}(x_1,\cdots, x_n)\,,\quad\mathbb{F}(x_1,\cdots, x_{n-1})\;\text{ with } \mathbb{F} \text{ finite of odd characteristic}
 \] have UUD\,?
\end{ques}

\begin{remark}\label{4.5}
Consider the group $k^*/k^{*2}$ of nonzero square classes in $k$. Any complete system of representatives for $k^*/k^{*2}$ is trivially a universality criterion set for quadratic forms over $k$. So we have a finite criterion set for forms over $k$ whenever this group of square classes is finite.

Notice that for some fields $k$ with  $k^*/k^{*2}$ infinite, it may still be possible to find such a finite criterion set. For example, if $u(k)=2$, then for any two elements $\alpha,\beta\in k^*$ not in the same square class, the set $S=\{\alpha,\beta\}$ is a universality criterion set for forms over $k$.  

More generally, if $\kappa$ is a field of characteristic $\neq 2$ that admits a finite universality criterion set $S_0$, then as we have said in the proof of Lemma\;\ref{3.4} (2), $S_0\cup S_0t$ is a finite universality criterion set for $\kappa(\!(t)\!)$. Therefore, if in Theorem\;\ref{thm-main} we take $k=\mathbb{C}(t)(\!(t_1)\!)\cdots(\!(t_n)\!)$ for some $n\in\mathbb{N}$, then $k^*/k^{*2}$ is infinite, but there still exists a finite criterion set $S$ for universal quadratic forms over $k$, and hence $S\cup Sx$ is a finite universality criterion set for definite quadratic forms over $k[x]$.
\end{remark}

\begin{remark}
Let $k$ be a non-archimedean local field of residue characteristic $p$, and let $S$ be a complete system of representatives for the square classes in $k^*$. We claim that $S$ is a minimal universality criterion set for quadratic forms over $k$.

Indeed, we can choose a ternary anisotropic form $\phi$ over $k$. For any $u\in S$, let $\varphi=(-u\det(\phi)).\phi$ be the scalar multiple of $\phi$ by $-u\det(\phi)$.
Then $\det(\varphi)=-u$, and it follows from \cite[63:21]{OMeara00} that $\varphi$ represents all the elements in $S\setminus\{u\}$ but does not represent $u$.

We leave it to the reader to check that the above statement can be generalized to higher local fields.
\end{remark}

\medskip

\noindent \emph{AI disclosure.} We selected and began the research topic of this paper in 2023, and in 2024 we obtained the results in \S\;2 and 3 of this paper. But later our research stalled because we wanted to use a Hasse principle approach which failed to work. In July 2026, the first named author asked ChatGPT-5.6 whether there was a way to handle the local field case. ChatGPT provided the idea of using ADC forms, which were studied in \cite{Cla12}. Using this idea, we generalized the local field case and ultimately proved the main theorem of this paper. This paper was written entirely by the authors, who take full responsibility for the content of the paper.

\medskip

\noindent \emph{Acknowledgements.}  We would like to thank Professors Ting-Yu Lee and Fei Xu for helpful discussions. The first named author was supported by a grant from the National Natural Science Foundation of China (no.\,12171223). The second named author benefited from the National Natural Science Foundation of China (no.\,12371063), and he thanks Professor Xiaowen Hu for kind support.

\addcontentsline{toc}{section}{\textbf{References}}

\bibliographystyle{alpha}

\bibliography{ref}

\begin{thebibliography}{BGVG14}

\bibitem[BGVG14]{BGVG14}
Karim~Johannes Becher, David Grimm, and Jan Van~Geel.
\newblock Sums of squares in algebraic function fields over a complete
  discretely valued field.
\newblock {\em Pacific J. Math.}, 267(2):257--276, 2014.

\bibitem[BH11]{BH11}
M.~Bhargava and J.~Hanke.
\newblock Universal quadratic forms and the 290-theorem.
\newblock Preprint, see
  \url{https://math.stanford.edu/\~vakil/files/290-Theorem-preprint.pdf}, 2011.

\bibitem[Bha00]{Bha00}
Manjul Bhargava.
\newblock On the {C}onway-{S}chneeberger fifteen theorem.
\newblock In {\em Quadratic forms and their applications ({D}ublin, 1999)},
  volume 272 of {\em Contemp. Math.}, pages 27--37. Amer. Math. Soc.,
  Providence, RI, 2000.

\bibitem[Bur06]{Bur06}
Jean~Edouard Bureau.
\newblock {\em Representation properties of definite lattices in function
  fields}.
\newblock ProQuest LLC, Ann Arbor, MI, 2006.
\newblock Thesis (Ph.D.)--Louisiana State University and Agricultural \&
  Mechanical College.

\bibitem[Cas26]{Cassady26JPAA}
Connor Cassady.
\newblock Universal quadratic forms over semi-global fields.
\newblock {\em J. Pure Appl. Algebra}, 230(7):Paper No. 108297, 32, 2026.

\bibitem[CD05]{CD05}
Wai~Kiu Chan and Joshua Daniels.
\newblock Definite regular quadratic forms over {${\Bbb F}_q[T]$}.
\newblock {\em Proc. Amer. Math. Soc.}, 133(11):3121--3131, 2005.

\bibitem[Cla12]{Cla12}
Pete~L. Clark.
\newblock Euclidean quadratic forms and {ADC} forms: {I}.
\newblock {\em Acta Arith.}, 154(2):137--159, 2012.

\bibitem[Con00]{Con00}
J.~H. Conway.
\newblock Universal quadratic forms and the fifteen theorem.
\newblock In {\em Quadratic forms and their applications ({D}ublin, 1999)},
  volume 272 of {\em Contemp. Math.}, pages 23--26. Amer. Math. Soc.,
  Providence, RI, 2000.

\bibitem[Djo76]{Djo76}
Dragomir\v{Z}. Djokovi\'c.
\newblock Hermitian matrices over polynomial rings.
\newblock {\em J. Algebra}, 43(2):359--374, 1976.

\bibitem[Ger03]{Ger03}
Larry~J. Gerstein.
\newblock Definite quadratic forms over {${\Bbb F}_q[x]$}.
\newblock {\em J. Algebra}, 268(1):252--263, 2003.

\bibitem[Ger04]{Ger04}
Larry~J. Gerstein.
\newblock On representation by quadratic {$\Bbb F_q[x]$}-lattices.
\newblock In {\em Algebraic and arithmetic theory of quadratic forms}, volume
  344 of {\em Contemp. Math.}, pages 129--134. Amer. Math. Soc., Providence,
  RI, 2004.

\bibitem[GVG92]{GesquiereVanGeel92}
Natacha Gesqui\`ere and Jan Van~Geel.
\newblock Note on universal quadratic forms.
\newblock {\em Bull. Soc. Math. Belg. S\'er. B}, 44(2):193--205, 1992.

\bibitem[HHK09]{HHK09}
David Harbater, Julia Hartmann, and Daniel Krashen.
\newblock Applications of patching to quadratic forms and central simple
  algebras.
\newblock {\em Invent. Math.}, 178(2):231--263, 2009.

\bibitem[KWX08]{KWX08}
Myung-Hwan Kim, Yuanhua Wang, and Fei Xu.
\newblock Universal quadratic forms over polynomial rings.
\newblock {\em J. Korean Math. Soc.}, 45(5):1311--1322, 2008.

\bibitem[Lam05]{Lam05}
T.~Y. Lam.
\newblock {\em Introduction to quadratic forms over fields}, volume~67 of {\em
  Graduate Studies in Mathematics}.
\newblock American Mathematical Society, Providence, RI, 2005.

\bibitem[Lee13]{Leep13}
David~B. Leep.
\newblock The {$u$}-invariant of {$p$}-adic function fields.
\newblock {\em J. Reine Angew. Math.}, 679:65--73, 2013.

\bibitem[O'M00]{OMeara00}
O.~Timothy O'Meara.
\newblock {\em Introduction to quadratic forms}.
\newblock Classics in Mathematics. Springer-Verlag, Berlin, 2000.
\newblock Reprint of the 1973 edition.

\bibitem[Sch85]{Schar}
Winfried Scharlau.
\newblock {\em Quadratic and {H}ermitian forms}, volume 270 of {\em Grundlehren
  der Mathematischen Wissenschaften [Fundamental Principles of Mathematical
  Sciences]}.
\newblock Springer-Verlag, Berlin, 1985.

\end{thebibliography}

\




















\end{document}